\documentclass[11pt]{amsart}
\usepackage{amscd,amsmath,amssymb}

\newcommand{\la}{\lambda}
\newcommand{\al}{\alpha}
\newcommand{\be}{\beta}

\newcommand{\f}{\varphi}

\newcommand{\norm}[1]{\Vert #1\Vert}

\newcommand{\Ll}{\mathcal{L}}
\newcommand{\vv}{\overrightarrow}

\newcommand{\Ker}{\mathop\mathrm{Ker\,}}

\newcommand{\CC}{\mathbb{C}}

\newcommand{\RR}{\mathbb{R}}
\newcommand{\ZZ}{\mathbb{Z}}

\numberwithin{equation}{section}

\newtheorem{te}{Theorem}[section]

\newtheorem{pr}{Proposition}[section]

\theoremstyle{definition}
\newtheorem{de}{Definition}[section]
\newtheorem{re}{Remark}[section]
\newtheorem{ex}{Example}[section]

\begin{document}

\title{Non-zero contact and  Sasakian reduction}
\thanks{Oana Dr\u agulete thanks the Swiss National Science Foundation
for partial support.
Liviu Ornea  was partially supported by ESI
 (Wien)
during August 2003, in the framework of the program "Momentum maps
and Poisson geometry"   and by Institut Bernoulli, EPFL, in July 2004, during the program ``Geometric mechanics and its applications''.}
\author{Oana Dr\u agulete}
\address{Section de math\'ematiques, EPFL, CH-1015 Lausanne, Switzerland
\newline\indent
and
\newline\indent
D{e}partement of Math{e}matics, \newline \indent University
"Politehnica" of
Bucharest, Romania}
\email{oana.dragulete@epfl.ch}

\author{Liviu Ornea}
\address{University of Bucarest, Faculty of Mathematics
\newline \indent
14 Academiei str.,
70109 Bucharest, Romania}
\email{Liviu.Ornea@imar.ro}


\begin{abstract}
We complete the reduction of Sasakian manifolds with the non-zero
case by showing that Willett's contact reduction is compatible with the
 Sasakian
structure. We then prove the compatibility of the non-zero Sasakian (in
particular, contact) reduction with the reduction of the K\"ahler (in
particular, symplectic) cone. We provide examples obtained by toric actions on
Sasakian spheres and make some comments concerning the curvature
of the quotients.

\end{abstract}

\maketitle

\section{Introduction}

\subsection{Sasakian manifolds}
We start by briefly recalling the notion of a Sasakian manifold,
sending to  \cite{bl} and \cite{bg_se} for more details and
examples.

\begin{de}
A {\it Sasakian} manifold is a $(2n+1)$-dimensional Riemannian manifold
 $(M,g)$
endowed with a unitary Killing vector field $\xi$ such that the
curvature tensor of $g$ satisfies the equation:
\begin{equation}\label{unu}
R(X,\xi)Y=\eta(Y)X-g(X,Y)\xi
\end{equation}
 where
$\eta$ is the metric
dual 1-form of $\xi$: $\eta(X)=g(\xi,X)$.
\end{de}

It can be seen that $\eta$ is a contact form (with Reeb field
$\xi$). Using the Killing property of $\xi$ and equation
\eqref{unu}, one defines an almost complex structure  on the
contact distribution $\Ker \eta$, by (the restriction of)
$\f=\nabla\xi$, where $\nabla$ is the Levi-Civita connection of
$g$.

The following formulae are then easily deduced:
\begin{equation}\label{doi}
\f \xi=0,\quad g(\f Y,\f Z)=g(Y,Z)-\eta(Y)\eta(Z).
\end{equation}

     The simplest compact example is the round  sphere
$S^{2n-1}\subset \CC^{n}$, with the metric induced by the flat one of
 $\CC^{n}$.
The characteristic  Killing vector field is $\xi_p=-i\vv{p}$, $i$ being
 the
imaginary unit. More general Sasakian structures on the sphere can be
 obtained
by deforming this standard structure as follows.
Let $\eta_A=\frac{1}{\sum a_j|z_j|^2}\eta_0$, for $0<a_1\leq
a_2\cdots\leq a_n$.  Its Reeb field is $R_A=\sum a_j(x_j\partial
y_j-y_j\partial x_j)$. Clearly, $\eta_0$ and $\eta_A$ underly the
same contact structure. Define the metric $g_A$ by the conditions:
\begin{itemize}
\item $g_A(X,Y)=\frac 12 d\eta_A(IX,Y)$ on the
  contact distribution (here $I$ is the standard complex structure of
 $\CC^n$);
\item $R_A$ is normal to the contact distribution and has unit length.
\end{itemize}
It can be seen  that $S^{2n-1}_A:=(S^{2n-1}, g_A)$ is a Sasakian
manifold (cf. \cite{ko}). It has recently been shown in \cite{ov} that each
compact
Sasakian manifold admits a CR-immersion in a $S_A^{2N-1}$.

Sasakian manifolds, especially the Sasakian-Einstein ones, seem to be
 more and
more important in physical theories (connected with the Maldacena
 conjecture).
Many new examples appeared lately, especially in the work of Ch.P.
 Boyer, K.
Galicki and their collaborators.

This growing importance of Sasakian structures was the first motivation
 for
extending in \cite{go} the contact (zero) reduction to this metric
 setting, by
showing that the contact reduction is compatible with the Sasakian data.

A good procedure for contact reduction away from zero was not
available
 when the paper \cite{go} was
written. We here complete the missing picture by showing that Willett's
 recently
defined non-zero
reduction introduced in \cite{wil} is compatible with the Sasakian data.

\bigskip

\noindent{\bf Acknowledgements.}~ The authors thank P. Gauduchon,
E. Lerman and T.S. Ratiu for their interest in this work and for
very useful suggestions and criticism.

\subsection{Contact reduction}

\subsubsection{Contact reduction at $0$ following \cite{al}, \cite{ge}}

Let $(M^{2n-1}, \eta)$ be an exact  contact manifold: this means that
 $\eta$ is
a contact form ($\eta
\wedge (d\eta)^n\ne 0$), hence its kernel is a contact structure on $M$.

 Let
 $R$ be the Reeb vector field, characterized by the conditions
$\eta (R) = 1$ and $d\eta (R, \cdot ) = 0$. The flow of the
(nowhere vanishing) Reeb vector
 field
preserves the contact form $\eta$.

Let $\Phi:G\times M\rightarrow M$ be an action by strong
 contactomorphisms of a
 (finite dimensional) Lie group on $M$: for any
$f\in G$, $f^*\eta=\eta$\footnote{If the action is proper or $G$ is
 compact,
 this is not
more restrictive than asking $G$ to preserve only the contact
structure: in the first case, one uses a Palais type argument, in
the second case an
 invariant
contact form can be found by averaging.}.
Such a  $G$--action by  strong
contactomorphisms on $(M, \eta)$ always admits an equivariant momentum
map $J:M\rightarrow \mathfrak{g}^*$ given by evaluating the contact
form on fundamental fields: $\langle J, \xi\rangle =
 \eta(\xi_M)$\footnote{Here
and in the sequel, for a $X\in \mathfrak{g}$, $X_M$ denotes the
 fundamental
field it induces on $M$.}. Note the main
difference towards the symplectic case: an action by contactomorphisms
 is automatically Hamiltonian.

It can be seen that $0\in \mathfrak{g}^*$ is  a regular value for
$J$ if and only if the fundamental fields induced by the action do
not vanish on the zero level set of $J$. In this case, the pull
back of the contact form to $J^{-1}(0)$ is basic. Let
$\pi_0:J^{-1}(0)\rightarrow J^{-1}(0)/G$ and
$\iota_0:J^{-1}(0)\hookrightarrow M$ be the canonical projection
(we shall always suppose that the considered actions are free and
proper, although these hypothesis can be relaxed to deal with the
category of orbifolds) and inclusion respectively. Albert's
reduction theorem assures the existence of a unique contact form
$\eta_0$ on $J^{-1}(0)/G$ such that $\pi_0^*\eta_0=\iota_0^*\eta$.
It can be seen that the contact structure of the quotient depends
only on the contact structure on $M$.

The Sasakian version of this result states (cf. \cite{go}) that if $M$ is Sasakian and
 $G$ acts
by isometric strong contactomorphisms, then the metric also projects to
 the
contact quotient and the whole structure is Sasakian.

\subsubsection{Contact reduction away from zero following \cite{wil}}
For $\mu\neq 0$, the restriction of the Reeb field is no more
basic on $J^{-1}(\mu)$ with respect to the action of $G_\mu$,
hence the above scheme does not apply. This situation was
corrected by Albert, but in an unsatisfactory way, see \cite{wil}
for examples. Willett's method, that
 we now
describe, is
more appropriate and was already used in \cite{dor} to extend the
 cotangent
reduction theorems in the contact context.  In the above setting,
for a  $\mu \in \mathfrak{g}^*$,
Willett calls the \emph{kernel group of $\mu$}, the connected Lie
subgroup $K_\mu$
of $G_\mu$  with Lie algebra $\mathfrak{k}_\mu = \ker\:
(\mu| _{\mathfrak{g}_\mu})$. One can see
that $\mathfrak{k}_\mu$ is an ideal in $\mathfrak{g}_\mu$, hence
 $K_\mu$ is a connected normal subgroup of $G_\mu$.
The contact quotient
of $M$ by $G$ at $\mu$ is defined by Willett as
$$
M_\mu : = J^{-1} (\RR_+ \mu)/K_\mu.
$$

\medskip
\noindent If $K_\mu$
acts freely and properly on $J^{-1} (\RR_+ \mu)$, then $J$ is
transversal to $\mathbb{R}_+ \mu$ and the pull back
of $\eta$ to  $J^{-1} (\RR_+ \mu)$ is basic relative to the
$K_\mu$--action on $J^{-1} (\RR_+\mu)$, thus inducing
a 1-form $\eta_\mu$ on the quotient $M_\mu$. If, in addition,
$\ker \mu + \mathfrak{g}_\mu = \mathfrak{g}$ then the form $\eta_\mu$
is also a contact form. It is characterized, as usual, by the identity
$\pi_\mu^* \eta_\mu = i_\mu^* \eta$, where
$\pi_\mu: J^{-1} (\RR_+ \mu) \to M_\mu$ is the canonical projection
and $i_\mu: J^{-1} (\RR_+ \mu) \hookrightarrow N$ is the canonical
 inclusion.

\begin{re} For $\mu=0$, Albert's and Willett's quotients
coincide.
\end{re}

In the next section we prove the compatibility of this procedure
with the metric context.

\section{Main results}
\subsection{The reduction theorem}
\begin{te}\label{1}
Let  $(M,g,\xi,\eta)$ be a $(2n-1)$ dimensional Sasakian manifold,
let $G$
be a Lie group of dimension $d$ acting on $M$ by strong contactomorphisms.
 Let $J:M\rightarrow\mathfrak{g}^*$
be the  momentum map associated to the action of $G$ and let   $\mu$ be
 an
{e}l{e}ment of the dual $\mathfrak{g}^*$.
We assume that:

$1.$ $\Ker\mu+\mathfrak{g}_\mu=\mathfrak{g}$.

$2.$ The action of $K_\mu$ on $J^{-1}(\RR_+\mu)$ is proper
and by isometries.

$3.$ $J$ is transverse to $\RR_+\mu$.

\noindent Then the contact quotient
$$
 M_\mu = J^{-1} (\RR_+\mu)/K_\mu
$$
is a  Sasakian manifold with respect to the projected metric and Reeb field.
\end{te}

\begin{proof}
We already know that the reduced space
 $M_\mu$ is a contact manifold (see \cite{wil}).
What is left to be proved is that the metric  $g$ and the  Reeb field
 $\xi$
 project
on $M_\mu$, the latter in a Killing field such that the curvature
 tensor
of the projected metric satisfies  formula \eqref{unu}.

From the transversality  condition satisfied by  the momentum map one
knows  that  $J^{-1}(\RR_+\mu)$ is an isometric Riemannian
submanifold of $M$ (which  induced
 metric we also denote by  $g$).
As the flow of the  Reeb field leaves invariant the level sets of the
 momentum
$J$, one derives that the restriction of  $\xi$ is still a unit Killing
 field on
 $J^{-1}(\RR_+\mu)$.

In order to establish the metric properties of
the canonical projection
$\pi_\mu:J^{-1}(\RR_+\mu)\rightarrow M_\mu$, we have to understand
the extrinsic geometry of the submanifold
$J^{-1}(\RR_+\mu)\subset M$. The first step is  to find a basis
in the normal bundle of
$J^{-1}(\RR_+\mu)$. To this end we look at  the
 direct sum  $\mathfrak{g}=\mathfrak{g}_{\mu}\oplus \mathfrak{m}$
 where $\mu\mid_{\mathfrak{m}}=0$ (such a decomposition
  exists, because
 $\Ker\mu+\mathfrak{g}_\mu=\mathfrak{g})$. Let
  $\mathfrak{m}_{M}=\{X_{M}\vert X\in \mathfrak{m}\}$ and recall that
(see \cite[Theorem 1]{wil}):
\begin{equation}
(T_{x}J^{-1}(\RR_+\mu)\cap \Ker \eta_{x})\oplus \RR \xi_x\oplus
{\mathfrak{m}}_{M}(x)=(T_{x}{\Phi}^{-1}(0)\cap \Ker \eta_x)\oplus
\RR \xi_x,
\end{equation}
for any  $x\in J^{-1}(\RR_+\mu)$, where  $\Phi$ is the  momentum
map associated to the action of $K_{\mu}$ on $M$.

Let now $\{X_{1},\ldots,X_{k}\}$ and $\{Y_{1},\ldots,Y_{m}\}$ be two
bases  in
$\mathfrak{k}_\mu$ and, respectively, $\mathfrak{m}$. Without loss of
generality, one may suppose that
the fundamental fields  $\{Y_{jM}\}_{j=1,m}$ form an  orthogonal basis
 of
$\mathfrak{m}_{M}$, $g$-orthogonal on $TJ^{-1}(\RR_+\mu)\cap \Ker \eta$
 and that
 $\{X_{iM}\}_{i=1,k}$ are  mutually orthogonal.

With these hypotheses, one derives that  $\{\varphi X_{iM},\varphi
  Y_{jM}\}$ are
 linearly ind{e}pendent in each $x\in M$  and
\[
g(\varphi Y_{jM},W)=g(\varphi X_{iM},W)=d\eta (W,X_{iM})=
-\langle dJ(W),X_{i}\rangle=\langle r\mu,X_{i}\rangle=0
\]
for any vector field  $W$  tangent to $J^{-1}(\RR_+\mu)$.
Therefore, for any
 $i,j$, the fields $\{ {\varphi X}_{iM},{\varphi Y}_{jM} \}$
belong to the normal bundle of $J^{-1}(\RR_+\mu)$. A simple counting of
 the
dimensions in the relation \eqref{1}, together with the fact that
 $\{{\varphi X}_{iM}\}$ is a basis in the normal bundle of
 $T\Phi^{-1}(0)$ (see the proof of  \cite[Theorem 3.1]{go}), imply that
$\{{\varphi X}_{iM},{\varphi Y}_{jM}\}$ is indeed a basis of the
normal bundle of $J^{-1}(\RR_+\mu)$\footnote{$\{{\varphi X}_{iM},{Y}_{jM}\}$ is also a basis for $T^\perp J^{-1}(\RR_+\mu)$ Our choice is only technically motivated.}.

Let $\nabla$, $\nabla^{M}$ be the Levi-Civita covariant
derivatives of $J^{-1}(\RR_+\mu)$ and $M$ respectively and let
$A_{i}$, $A_j$ be the Weingarten operators associated to the
unitary
 normal  sections
 ${\varphi X}_{iM}/\Vert X_{iM}\Vert$, $1\le i\le
k$, $\varphi Y_{jM}/\Vert Y_{jM}\Vert$, $1\le j\le m$.
By applying the  Weingarten formula and the relation
\eqref{unu}, one obtains, for any  $X,Y,Z$ tangent to
 $J^{-1}(\RR_+\mu)$:
\begin{equation*}
\begin{split}
g(A_{i}Y,Z)&={\Vert X_{iM}\Vert}^{-1}\{ g(X_{iM},Y)\eta (Z)-g(\varphi
\nabla^{M}_{Y}X_{iM},Z)\},\\
g(A_{j}Y,Z)&={\Vert Y_{jM}\Vert}^{-1}\{ g(Y_{jM},Y)\eta (Z)-g(\varphi
\nabla^{M}_{Y}Y_{jM},Z)\}.
\end{split}
\end{equation*}

As  $K_\mu$ acts by strong contact isometries, the
metric $g$  projects on a  metric $g^{M_\mu}$ on $M_\mu$
with respect to which
the canonical projection $\pi_\mu$ becomes a Riemannian  submersion.
We now show that the vertical distribution $\mathcal{V}$ is
locally generated by the vector fields
 $\{X_{iM}\}$. We have indeed:
\[
T_{x}\pi_{\mu}(X_{iM}(x))=
T_{x}\pi_{\mu}\left(\dot{c}(0) \right) =
\dot{(\pi_{\mu}\circ c)}(0)
\]
where $c(t)=\Phi(\exp tX_{iM},x)$.

But $(\pi_{\mu}\circ c)(t)=\pi_{\mu}(x)$ for any  $t$
and then
$$T_{x}\pi_{\mu}(X_{iM}(x))=0\quad  \text{for any  $x\in
 J^{-1}(\RR_+\mu)$}.$$
This
proves that
$\{X_{iM}\}_{1\le i\le k}\subset \mathcal{V}_{x}$ and,
as $\dim \mathcal{V}_{x}=k$, it implies that
$\{X_{iM}\}$ generate $\mathcal {V}$.

The formulae $\Ll_{X_{iM}}\xi=0$ for  $i=1,\ldots,k$
prove that $\xi$ is a projectable vector field and its projection
$\zeta$ is a unit Killing field
on the reduced space $M_{\mu}$.

Let  $X,Y,Z$ be vector fields  orthogonal to $\zeta$. Using  O'Neill's
 formulae
(see \cite[(9.28f)]{bes})
we derive:
\begin{equation*}
\begin{split}
g^{M_{\mu}}(R^{M_\mu}(X,\zeta)Y,Z)&=
g(R(X^{h},\xi)Y^{h},Z^{h})+2g(A(X^{h},\xi),A(Y^{h},Z^
{ h})\\
&-g(A(\xi,Y^{h}),A(X^{h},Z^{h}))+g(A(X^{h},Y^{h}),A(\xi,Z^{h})),
\end{split}
\end{equation*}
where $X^{h}$ denotes the horizontal lift of the vector field $X$,
$A$ is O'Neill's (1,2) tensor field given by the relation:
$A(Z^h,X^h)=\text{\emph{vertical part of}}\, \nabla^M_{Z^h}X^h$ and $R$ the curvature tensor of the connection $\nabla$ on $J^{-1}(\RR_+\mu)$.
On the other hand:
\begin{equation*}
\begin{split}
g(\nabla_{Z^{h}}\xi,X_{iM})&=g(\varphi
 Z^{h},X_{iM})=d\eta(X_{iM},Z^{h})=\langle
dJ(Z^{h}),X_{iM}\rangle\\
&=r\langle\mu,X_{iM}\rangle=0,
\end{split}
\end{equation*}
and hence:
\[
R^{M_{\mu}}(X,\zeta)Y=R(X^{h},\xi)Y^{h}.
\]
This completes the proof.

\end{proof}

\begin{re}
Under the  hypothesis of the above theorem, the dimension of the
 reduced
space is  $2n-d-m-k+1$.
\end{re}

\subsection{Compatibility with the K\"ahler reduction}
We now analyze the
compatibility of the non-zero Sasakian reduction with  K\"ahler
reduction using the cone construction. In particular, we  obtain a
relation between non-zero contact reduction and symplectic  reduction.

Let $\mathcal{C}(M)=M\times \RR_+$ be the cone over $M$ endowed with
 the K\"ahler metric $r^2g+dr^2$.
The action of $G$ on $M$ lifts to an action on $\mathcal{C}(M)$ by
 holomorphic isometries which commute with the translations along the
 generators (see \cite{go} e.g.).
Similarly, the action of $K_\mu$ lifts to the cone, the lifted action
 being the restriction of the above.

Let $J_s:\mathcal{C}(M)\rightarrow \mathfrak{g}^*$, resp.
 $\Phi_s:\mathcal{C}(M)\rightarrow \mathfrak{k_\mu}^*$ be the symplectic
 momentum map associated to the $G$-action, resp. $K_\mu$-action on the
 cone.  The differentials of these two momentum maps are related by the
transpose $\iota^t$ of the natural  inclusion
$\iota:\mathfrak{k_\mu}\hookrightarrow \mathfrak{g}$, namely
$T\Phi_s=\iota^t\circ TJ_s.$

We now embed $M$ in the cone as $M\times \{1\}$ and observe that the
contact momentum maps are the restrictions of the symplectic ones:
$J=J_s|_{M\times\{1\}}$, resp. $\Phi=\Phi_s|_{M\times\{1\}}$. Clearly
$J$ and $\Phi$ are the contact momentum maps associated to the $G$,
resp. $K_\mu$-action on $M\times\{1\}$. Moreover, we have
\begin{equation}\label{iii}
\Phi=\iota^t\circ J.
\end{equation}

On the other hand, we recall (see \cite{go}) that the reduced space at $0$
 of the K\"ahler cone is the K\"ahler cone of the Sasakian reduced
space
 at $0$:
$$\Phi_s^{-1}(0)/K_\mu=\mathcal{C}(\Phi^{-1}(0)/K_\mu).$$

We are now prepared to prove:
\begin{te}\label{comp}
Let  $(M,g,\xi,\eta)$ be a Sasakian manifold,
let $G$
be a Lie group acting on $M$ by strong contactomorphisms
and $\mu$ an element of $\mathfrak{g}^*$. Suppose that:
\begin{itemize}
\item $0$ is a regular value for $J_s$.
\item $\Ker \mu+\mathfrak{g}_\mu=\mathfrak{g}$.
\item $K_\mu$ acts properly and by isometries on $\Phi^{-1}(0)$.
\item $J$ is transverse to $\RR_+\mu$ and to $\RR_-\mu$.
\end{itemize}
Then the cone over the Sasakian quotient of $M$ at $0$ with respect to
 the $K_\mu$ action is the disjoint union of the K\"ahler cones over the
 Sasakian quotients $M_\mu$ and $M_-\mu$ and the cone over the
 co-isotropic submanifold $J^{-1}(0)/K_\mu$:
$$(\mathcal{C}(M))_0=\mathcal{C}(\Phi^{-1}(0)/K_\mu)=
\mathcal{C}(M_\mu)\cup \mathcal{C}(J^{-1}(0)/K_\mu)\cup
\mathcal{C}(M_-\mu).$$
\end{te}
\begin{proof}
Since $\iota^t$ is surjective, from \eqref{iii} and
from $0$ being a regular value for $J_s$ (and hence also for $J$) it
follows that ${0}$ is a regular value for $\Phi_s$ and hence for
$\Phi$.

As $\Phi^{-1}(0)=J^{-1}(\RR\mu)$ (cf. \cite[proof of Theorem 2]{wil})
 and $K_\mu$ acts on
$J^{-1}(0)$  and on $J^{-1}(\RR_+\mu)$, we have the partition
$$\Phi^{-1}(0)/K_\mu=J^{-1}(\RR_+\mu)/K_\mu\cup J^{-1}(0)/K_\mu\cup
J^{-1}(\RR_-\mu)/K_\mu.$$
We note that $(M,g,-\xi)$ is also a Sasakian manifold on which $G$ acts
 by Sasakian automorphisms and the associated momentum map is $-J$.
 Then, if  the quotient $M_{-\mu}:=J^{-1}(\RR_-\mu)/K_\mu$ exists (or,
equivalently, the quotient of $(M,g,\xi)$ at $-\mu$), it
 will be a Sasakian manifold according to our previous theorem. But note
 that two of the hypothesis of the theorem are not automatically
 satisfied in both cases: if $K_\mu$ acts properly on $J^{-1}(\RR_+\mu)$
 it does not necessarily act properly on $J^{-1}(\RR_-\mu)$ and
 similarly for the transversality condition.

Now $J^{-1}(0)/K_\mu$ is a manifold on which the one-form $\eta$ is
projected. However, it is no more a contact form, nor, in general, a
symplectic one. Indeed, using Albert \cite[Propositions
1,2]{al}, $J^{-1}(0)/K_\mu$ is contact or symplectic if and only if
$$T_x(K_\mu\cdot x)=\Ker
(T\eta_x|_{T_xJ^{-1}(0)\cap \Ker \eta_x}).$$
But, in general, one has $\Ker
(T\eta_x|_{T_xJ^{-1}(0)\cap \Ker \eta_x})=T_x(G\cdot x)$.
However, this implies that $J^{-1}(0)/K_\mu$ is a co-isotropic
submanifold with respect to the contact form of $\Phi^{-1}(0)/K_\mu$.
 \end{proof}

\begin{re}
Forgetting the metric and \emph{mutatis mutandis}, the result of Theorem \ref{comp} remains valid for
contact manifolds.
\end{re}

\section{Examples: actions of tori on spheres.}

In Willett's reduction scheme, the smallest dimension of $G$ which
produces non-trivial examples is $2$. We here present some
complete computations for various actions of $G=T^2$ on $M=S^7$
with the standard Sasakian structure given by the contact form
$\eta= \sum  (x_{j}dy_{j}-y_{j}dx_{j})$.  When possible, we
briefly discuss also the reduction at zero with the same group and the cone
construction (the notations for the momentum maps will be the ones used in the
previous section).
Generalizations to $S^{2n-1}$ are also indicated.

Note that our examples show the dependence of the dimension of the
quotient on the choice of $\mu$.

\begin{ex}
Let first $T^2$ act on $S^7$  by
 $$((e^{it_{0}},e^{it_{1}}),(z_{0},...,z_{3}))\mapsto
(e^{it_{0}}z_{0}, e^{it_{0}}z_{1}, e^{it_{1}}z_{2},
 e^{it_{1}}z_{3}).$$
 Since $G$ is
 commutative, $\mathfrak{g}_{\mu}=\mathfrak{g}=\mathbb{R}^{2}$.

For any $(r_{1},r_{2})\in \mathfrak{g}$ the associated infinitesimal
 generator is given by
\begin{equation*}
\begin{split}
(r_{1},r_{2})_{S^{7}}(z)&=
r_{1}(-y_{0}\partial_{x_{0}}+x_{0}\partial_{y_{0}})+r_{1}(-y_{1}\partial
_{x_{1}}+x_{1}\partial_{y_{1}})\\
&+r_{2}(-y_{2}\partial_{x_{2}}+x_{2}\partial_{y_{2}})+
r_{2}(-y_{3}\partial_{x_{3}}+x_{3}\partial_{y_{3}})
\end{split}
\end{equation*}
and the momentum map
$J:S^7\rightarrow {(\RR^2)}^*$ reads $J(z)=\langle(\vert  z_{0}\vert^2
+\vert z_{1}\vert^2 ,\vert z_{2}\vert^2 +\vert
 z_{3}\vert^2),\cdot\rangle$.

Let $\mu:\RR^2\rightarrow \RR$, $\mu=\langle v,\cdot\rangle$, $v\in
 \RR^2\setminus\{0\}$ fixed. Then:
\begin{equation*}
J^{-1}(\RR_+\mu)=
\begin{cases}
S^3(\sqrt{\frac{v_1}{v_1+v_2}})\times S^3(\sqrt{\frac{v_2}{v_1+v_2}}),
 &\text{if $v_1,v_2>0$}\\
S^3(\sqrt{\frac{v_1}{v_1+v_2}}),
 &\text{if $v_1>0,v_2=0$}\\
S^3(\sqrt{\frac{v_2}{v_1+v_2}}),
 &\text{if $v_1=0,v_2>0$}\\
\end{cases}
\end{equation*}

For $v=(1,0)$ $J^{-1}(\RR_+\mu)=S^3$, $\Ker
\mu=\mathfrak{k}_\mu=\{0\}\times\RR$,  $K_\mu=\{e\}\times S^1$. The action
of $K_\mu$ on  $J^{-1}(\RR_+\mu)$ is  trivial and hence $M_\mu=S^3$.
In this case $0$ is
not a regular value of $\Phi$-the momentum map associated to the $K_{\mu}$
action but, nevertheless, $\Phi^{-1}(0)$ is a submanifold of $S^7$ and hence
the reduced space at zero, $\Phi^{-1}(0)/K_{\mu}$ is a Sasaki manifold. As
$\Phi^{-1}(0)=S^3$ and  $\mathcal{C}(S^n)=\RR^{n+1}\setminus\{0\}$,
we obtain that $(\mathcal{C}(S^7))_0=\RR^4\setminus \{0\}$. Note that for
this choice of $\mu$
reducing and taking the cone are commuting operations exactly as in the
zero case.

For $v=(1,1)$ we obtain: $J^{-1}(\RR_+\mu)=S^3(\frac{1}{\sqrt2})\times
S^3(\frac{1}{\sqrt2})$,  $\mathfrak{k}_\mu=\{(-x,x)\vert x\in \RR\}$,
$K_\mu= \{(e^{-it},e^{it})\vert e^{it}\in S^1\}$. The action of $K_\mu$
on  $J^{-1}(\RR_+\mu)$ is given by
$$((e^{-it},e^{it}),z)\mapsto
(e^{-it}z_0,e^{-it}z_1,e^{it}z_2,e^{it}z_3),$$
thus $M_\mu=S^2\times S^3$.

We can generalize this example for $M=S^{2n+1}$ by considering the
action
$$((e^{it_0},e^{it_1}),z)=(e^{it_0}z_0,e^{it_0}z_1,e^{it_1}z_2,...,e^{it
_1}z_n).$$
Now the momentum map is $J(z)=\langle(\vert z_{0}\vert^2 +\vert
z_{1}\vert^2,\sum{\vert z_k\vert^2}),\cdot\rangle$.
For $\mu$ as above, we have:
\begin{equation*}
J^{-1}(\RR_+\mu)=
\begin{cases}
S^3(\sqrt{\frac{v_1}{v_1+v_2}})\times
S^{2n-3}(\sqrt{\frac{v_2}{v_1+v_2}}),  &\text{if $v_1,v_2>0$}\\
S^3(\sqrt{\frac{v_1}{v_1+v_2}}),
 &\text{if $v_1>0,v_2=0$}\\
S^{2n-3}(\sqrt{\frac{v_2}{v_1+v_2}}),
 &\text{if $v_1=0,v_2>0$}\\
\end{cases}
\end{equation*}

 For the same particular choices of $\mu$ as above, we obtain as reduced
spaces respectively $S^3$,  $S^{2n-3}$ or
$S^3\times\mathbb{C}P^{n-2}$.

\end{ex}

\begin{ex}
Let now the action be given by
$$((e^{it_{0}},e^{it_{1}}),z)\mapsto
(e^{-it_{0}}z_{0}, e^{it_{0}}z_{1}, e^{it_{1}}z_{2}, e^{it_{1}}z_{3}).$$
The infinitesimal generator of the action will be
\begin{equation*}
\begin{split}
(r_1,r_2)_{S^7}(z)&=r_{1}(y_{0}\partial_{x_{0}}-x_{0}\partial_{y_{0}}) +
r_{1}(-y_{1}\partial_{x_{1}}+x_{1}\partial_{y_{1}})\\
&+r_{2}(-y_{2}\partial_{x_{2}}+x_{2}\partial_{y_{2}})+
r_{2}(-y_{3}\partial_{x_{3}}+x_{3}\partial_{y_{3}}).
\end{split}
\end{equation*}
The momentum map is
$J(z)=\langle(\vert z_1\vert^2-\vert z_0\vert^2,\vert z_2\vert^2+\vert
 z_3\vert^2),\cdot\rangle$ and
\begin{equation}\label{es}
J^{-1}(\RR_+\mu)=\left\{z\in S^7\,\mid\, \exists
s>0\,\,\text{such that}
\begin{cases}
\vert z_1\vert^2-\vert z_0\vert^2-sv_1=0,\\
\vert z_2\vert^2+\vert z_3\vert^2-sv_2=0.\\
\end{cases}
\right\}
\end{equation}
For $v=(1,0)$ we obtain
$$J^{-1}(\RR_+\mu)=\{z\in S^7\vert
 z_2=z_3=0,\vert z_1\vert>\vert z_0\vert\}=
S^3\setminus\{\vert z_1\vert\leq\vert z_0\vert\}.$$
The action of
 $K_\mu=\{e\}\times S^1$ on $J^{-1}(\RR_+\mu)$
is trivial, thus $M_\mu=S^3\setminus\{\vert z_1\vert\leq\vert
 z_0\vert\}$, an open submanifold of $S^3$.
For $v=(1,1)$, solving for $s$ the equations in \eqref{es} gives
$s \in (0, \frac{1}{2}]$. Hence:
 \begin{equation*}
\begin{split}
J^{-1}(\RR_+\mu)&\simeq
\left(S^1(\frac{1}{\sqrt2})\times S^5(\frac{1}{\sqrt 2})\right)\setminus
\left\{z\in S^7\mid \vert z_0\vert^2=\frac12\right\}\\
&\simeq
 S^1(\frac{1}{\sqrt2})\times(S^5(\frac{1}{\sqrt 2})\setminus
 S^1(\frac{1}{\sqrt2}))
\end{split}
\end{equation*}
an open submanifold of the product of spheres.

 The action of $K_\mu$ on
 $J^{-1}(\RR_+\mu)$ is given by
$$((e^{-it},e^{it}),z)\mapsto
(e^{it}z_{0}, e^{-it}z_{1}, e^{it}z_{2}, e^{it}z_{3}).$$
Let $A$ denote the set $\left\{z\in S^7(\sqrt2)\mid 0<\vert
z_2\vert^2+\vert  z_3\mid^2\le1\right\}$. Obviously, the above action of
$K_\mu$ can be understood on the whole $\CC^4$ and, as such, restricts
to an action on $A$. Then  $M_\mu$ is diffeomorphic  with $(S^1\times S^5)\cap
A/K_\mu$.
To identify the quotient, let $g:\left(S^1\times S^5\right)\cap A
\rightarrow \left(S^1\times S^5\right)\cap A$ be given by
$$(z_0,z_1,z_2,z_3)\mapsto (z_0,z_1^{-1},z_2,z_3).$$
$g$ induces a map from
$\left(\left(S^1\times S^5\right)\cap A\right) S^1$ (with  respect to
the diagonal action of $S^1$)  to
$\left(\left(S^1\times S^5\right)\cap A\right)/K_\mu$. The map
$$(z_0,...,z_3)\mapsto (\bar z_1z_0,z_1,\bar z_1z_2,\bar z_1z_3)$$ is a
diffeomorphism of $\left(S^1\times S^5\right)\cap A$ equivariant  with
respect to the diagonal action of $S^1$ and the action of $S^1$ on the
first  factor. Hence $M_\mu$
is diffeomorphic to $S^5(\frac{1}{\sqrt 2})\setminus
\mathrm{pr}\left\{ z\in  S^7\mid \vert z_0\vert^2
=\frac 12\right\}\simeq S^5(\frac1{\sqrt2})\setminus S^1(\frac
 1{\sqrt2})$,
where $\mathrm{pr}:\CC^4\rightarrow\CC^3$,
$\mathrm{pr}(z_0,\ldots, z_3)=(z_0, z_2,z_3)$.

If we change the action on $z_0$ with $e^{-ikt}z_0$, the reduced space
 will be the  above one quotiented by $\ZZ^k$ (see also
\cite[Example 4.2]{go}).

\smallskip

 \end{ex}

\begin{ex}
Let us take this time:
$$((e^{it_0},e^{it_1}),z)\mapsto
(e^{it_{0}}z_{0},e^{it_{1}}z_{1},e^{it_{ 1 } }
z_{2},e^{it_{1}}z_{3}),$$
whose infinitesimal generator is
$$(r_1,r_2)_{S^7}(z)=r_1(-y_0\partial x_0+x_0\partial
y_0)+r_2\sum_{j=1}^3  (-y_j\partial x_j+x_j\partial y_j).$$
The momentum map is:
$$J(z)=\langle(\vert z_0\vert^2,\vert z_1\vert^2+\vert z_2\vert^2+\vert
 z_3\vert^2),\cdot\rangle.$$
For $J^{-1}(\RR_+\mu)$ we obtain the following possibilities:
\begin{equation}
J^{-1}(\RR_+\mu)=
\begin{cases}
S^1(\sqrt{\frac{v_1}{v_1+v_2}})\times S^5(\sqrt{\frac{v_2}{v_1+v_2}}),
 &\text{if $v_1,v_2>0$}\\
S^5(\sqrt{\frac{v_2}{v_1+v_2}}),
 &\text{if $v_1=0,v_2>0$}\\
S^1(\sqrt{\frac{v_1}{v_1+v_2}}),
 &\text{if $v_2=0,v_1>0$}\\
\end{cases}
\end{equation}
In particular, for $v=(1,0)$, $M_\mu=S^1$, for
$v=(0,1)$, $M_\mu=S^5$ and for $v=(1,1)$  one obtains the same quotient
as in the preceding example.
\end{ex}

\begin{ex}
Considering the weighted action of $T^2$ on $S^7$ given this time by
$$((e^{it_0},e^{it_1}),z)\mapsto
(e^{it_{0}\la_0}z_0,e^{it_{1}\la_1}z_1, z_2, z_3),$$ one obtains
the momentum map
$$J(z)=\langle(\la_0\vert z_0\vert^2,\la_1\vert
 z_1\vert^2),\cdot\rangle.$$
For $v=(0,1)$  and $\la_1$ strictly positive, the reduced space is
$S^5 \setminus S^3$ if $\la_0 \ne 0$ and $S^7 \setminus S^5$ if
 $\la_0=0$.

The cone construction is verified in this case. Indeed, $J^{-1}(0)=S^3$ and
$$(\mathcal{C}(S^7))_0\simeq \mathcal{C}(S^5) = \mathcal{C}(S^3)\cup
\mathcal{C}(S^5\setminus S^3).$$

If $v=(1,1)$
and $\la_0,\la_1>0$,
\begin{equation}
J^{-1}(\RR_+\mu)=\left\{z\in S^7\mid \vert
 z_1\vert=\sqrt{\frac{\la_0}{\la_1}}\vert z_0 \vert, z_0\ne 0\right\}=
S^7\cap\left(\CC^\ast\times A\right)
\end{equation}
where $A$ is the ellipsoid of equation
$$\vert z_1\vert^2(1+\frac{\la_1}{\la_0})+\vert z_2\vert^2+\vert
 z_3\vert^2=1.$$
The action of $K_\mu$ on $J^{-1}(\RR_+\mu)$ is given by
$$((e^{-it},e^{it}),z)\mapsto(e^{-it\la_0}z_0,e^{it\la_1}z_1,z_2,z_3)$$
and the reduced space
$$M_\mu=\bigcup_{(z_2,z_3)\in
 \mathrm{pr}(J^{-1}(\RR_+\mu))}S^1(\be^{-\la_0}\al^{\la_1})\times\left\{(z_2,z_3)
\right\}$$
where $\mathrm{pr}:\CC^4\rightarrow\CC^2$,
$\mathrm{pr}(z_0,\ldots, z_3)=(z_2,z_3)$, $\be=\sqrt{\frac{\la_0(1-\vert
z_2\vert^2-\vert z_3\vert^2)}{\la_0+\la_1}}$ and
$\al=\sqrt{\frac{\la_1(1-\vert z_2\vert^2-\vert
 z_3\vert^2)}{\la_0+\la_1}}$.

\medskip

If $[z]= [z']$ in the reduced space then $z_2=z'_2$ and
$z_3=z'_3$. So let $(z_2,z_3)$ be fixed in $pr(J^{-1}(\RR_+\mu))$.
$z\in J^{-1}(\RR_+\mu)$ and $\mathrm{pr} (z)=(z_2,z_3)$ imply $\vert
z_0\vert=\al$ and $\vert z_1\vert=\be $. The action of $K_\mu$ on
$J^{-1}(\RR_+\mu)$ is in fact the diagonal action of $S^1$ on the
first two coordinates. Let $f:(S^1(\al)\times S^1(\be)\times
\{(z_2,z_3)\})/S^1\rightarrow S^1(\al^{\la_1}\be^{-\la_0})$ be the
map given by
$$[z]\mapsto z_0^{\la_1}z_1^{\la_0}.$$
One can easily check that $f$ is a diffeomorphism.

\medskip

In the previous examples, the Reeb flow on the reduced space is the restriction
of the canonical one of the standard sphere. In this latter case, we obtain a
non-standard Reeb flow.

We now write the flow of the Reeb field of the reduced contact form on  $M_\mu$
(for $v=(1,1)$). Let $r(t)=\left(\begin{smallmatrix} \cos t &-\sin t\\ \sin
t&\cos t\end{smallmatrix}\right)$, $Z=(z_0^2, z_0^3)^t$. Then the flow is
written as
$$\f_t=(Ae^{i(a+bt)}, R(t)Z),$$
where
\begin{equation*}
\begin{split}
A&=\norm{z_0^0}^{\la_1}\norm{z_0^1}^{\la_0},\\
a&=\la_1v_0+\la_0v_1, \, \text{with} \, v_0=\mathrm{arg}(z_0^0),
v_1=\mathrm{arg}(z_0^1),\\
b&=\la_1+\la_0,\\
R(t)&=\mathrm{diag}(r(t), r(t)).
\end{split}
\end{equation*}

\end{ex}

\section{The sectional curvature of the quotient}

\subsection{Contact CR submanifolds}
In order to evaluate the sectional curvature of the Sasakian reduced
 space, both
at 0 and away from 0, it will be convenient to place ourselves in
a slightly more general situation. We first recall the following
definition (see e.g. \cite{be}):
\begin{de}
Let $(M, g_M, \xi)$ be a Sasakian manifold. An isometric
submanifold $N$ is called \emph{contact CR} or
\emph{semi-invariant} if it admits two mutually orthogonal
distributions $D$ and $D^\perp$,  such that:
\begin{enumerate}
\item $TN$  decomposes orthogonally as: $TN= D\oplus D^\perp\oplus
\langle\xi\rangle$ and \item $\f D =D,\quad \f D^\perp \subseteq
T^\perp N$.
\end{enumerate}
\end{de}
We see that, in general, the normal bundle of the submanifold also
 splits into
two orthogonal distributions: $\f D^\perp$ and its orthogonal complement
that we
denote by $\nu$ and which is invariant at the action of $\f$. We
then have:
$$TM_{|N}=D\oplus D^\perp\oplus \langle\xi\rangle\oplus \f D^\perp\oplus
\nu.$$
For a vector field $V$ normal to $N$ we shall denote $\bar V$,
respectively
 $\tilde V$
its component in $\f D^\perp$, respectively  in $\nu$.

Such submanifolds have been extensively studied in the last thirty years.

Obviously,  very  natural examples are the level sets of Sasakian
 momentum maps.
To better mimic our situation, we moreover make the following:\\[1mm]
\noindent {\bf Assumption.} There exists a Riemannian submersion
$\pi:N\longrightarrow P$ over a Sasakian manifold $(P,g_P,\zeta)$ such
 that:
 \begin{enumerate}
\item $D\oplus \langle\xi\rangle$ represents the horizontal
distribution of the submersion; (and hence $D^\perp$ represents
the vertical distribution of the submersion); \item The two Reeb
fields are $\pi$-related: $\xi$ is basic and projects over
$\zeta$.
\end{enumerate}

This situation was already considered by Papaghiuc in \cite{pap}, on the model
of Kobayashi's paper \cite{kob} where the similar setting was discussed in
K\"ahlerian context.

\smallskip

Let  $\phi:=\nabla^P\zeta$ and observe that in our assumption we have
 $(\phi
X)^h=\f X^h$.

We want to relate the sectional curvature of planes
 generated by orthonormal
pairs $\{X,\phi X\}$, respectively $\{X^h,\varphi X^h\}$. This is
 usually known
as $\f $-sectional curvature, the analogue in Sasakian geometry of
 holomorphic sectional curvature; it completely determines the curvature
 tensor,
 cf. \cite{bl}, so it is worth having information about it.

We first apply (as in the proof of  Theorem \ref{1}) O'Neill's formula
 to
relate the curvatures of $N$ and $P$. For $X$ tangent to $P$ and
 orthogonal to
$\zeta$, (this is not restrictive, as the planes passing through the
 Reeb field
have sectional curvature $1$ on a Sasakian manifold), using the
anti-symmetry of the tensor
 $A$,
we obtain:
\begin{equation}\label{onil}
R^N(X^h,\f X^h, X^h,\f X^h)-R^P(X,\phi X, X,\phi X)=-3 \Vert A(X^h,\f
X^h)\Vert_N^2,
\end{equation}
where the sub-index refers to the norm with respect to $g^N$.

The next step is to apply the Gauss equation to the Riemannian
submanifold $N$ of
 $M$:
\begin{equation}\label{gau}
\begin{split}
R^M(X^h,\f X^h, &X^h,\f X^h)-R^N(X^h,\f X^h, X^h,\f X^h)\\
&=\Vert h(X^h,\f X^h)\Vert_M^2-g_M(h(X^h,X^h),h(\f X^h, \f X^h)).
\end{split}
\end{equation}
We now need to relate the tensors $A$ and $h$. To this end, we write
 $hE$, respectively
$vE$ for the horizontal, respectively  vertical part of a tangent
(to $N$) vector field $E$ and we first decompose
$$\nabla^M_{X^h}(\f Y^h)=h\nabla^M_{X^h}(\f Y^h)+A(X^h,\f Y^h)+h(X^h,\f
 Y^h).$$
Then we use  the formulae $(\nabla^M_E\f)F=\eta(F)E-g_M(E,F)\xi$ (see
\cite{bl}) and $(\nabla^M_E\f)F=\nabla^M_E(\f F)-\f \nabla^M_EF$ to
 express
$\nabla^M_{X^h}(\f Y^h)$. Finally,  equaling the
tangent and normal parts in the equation we obtain this way, we arrive
at the following relations:
\begin{equation}\label{rel1}
\begin{split}
A(X^h,\f Y^h)&=v\f h(X^h,Y^h),\\
h(X^h, \f Y^h)&=\f A(X^h,Y^h)+\f\widetilde{h(X^h,Y^h)}.
\end{split}
\end{equation}
Note that if $\f D^\perp=T^\perp N$ (\emph{i.e.} $\nu=\{0\}$), and
this is the case when $N$ is the $zero$ level set of a Sasakian momentum map, the above
relations simplify to:
\begin{equation}\label{rel2}
\begin{split}
A(X^h,\f Y^h)&=\f h(X^h,Y^h),\\
h(X^h, \f Y^h)&=\f A(X^h,Y^h).
\end{split}
\end{equation}
In the general case, from \eqref{rel1} we easily derive:
$$h(\f X^h,\f Y^h)=\overline{h(X^h,Y^h)}-\widetilde{h(X^h,Y^h)},$$
and hence
\begin{equation}\label{norm1}
g_M(h(\f X^h,\f Y^h), h(X^h,Y^h))=\Vert \overline{h(X^h,Y^h)}\Vert_M^2-
\Vert \widetilde{h(X^h,Y^h)}\Vert^2_M.
\end{equation}
From equation \eqref{doi} it follows that on the orthogonal
complement of $\xi$, the tensor $\f$ acts like
an isometry. Therefore, using again \eqref{rel1}, we derive:
\begin{equation}\label{norm2}
\begin{split}
\Vert {h(X^h,\f Y^h)}\Vert_M^2&=\Vert A(X^h,Y^h)\Vert^2_M+
\Vert \widetilde{h(X^h,Y^h)}\Vert^2_M,\\
\Vert A(X^h,\f Y^h)\Vert^2_M&=\Vert \overline{h(X^h,Y^h)}\Vert_M^2.
\end{split}
\end{equation}
Let us denote $K_\phi^P(X)$, respectively  $K_\f^M(X^h)$ the
sectional curvature of the plane  $\{X,\phi X\}$, respectively
$\{X^h, \f X^h\}$. Adding equations \eqref{onil},
 \eqref{gau} and using \eqref{norm1},
\eqref{norm2}, we finally obtain (taking again into account the
 anti-symmetry of
$A$):
\begin{equation}\label{final}
K_\phi^P(X)=K_\f^M(X^h)+4\Vert
{\overline{h(X^h,X^h)}}\Vert_M^2-2\Vert\widetilde{h(X^h,X^h)}\Vert_M^2.
\end{equation}

\subsection{The curvature of the quotient}
In general, from equation \eqref{final} one hopes to deduce the
positivity of the $\f$-sectional curvature of the quotient. This
depends on the extrinsic geometry of the level set, which is a
data additional to the reduction scheme: the second fundamental
form of the level set cannot be entirely expressed in terms of the
action. But in some particular cases, one is able to derive a
conclusion.

Obviously the simplest situation occurs when $J^{-1}(\RR_+\mu)$ is totally geodesic in $M$: then the $\f$-sectional curvatures of $M$ and $M_\mu$ are equal. But even our examples show that this is not always the case. In fact, one is only interested in the vanishing of $h(X^h, Y^h)$, which, by the first equation in \eqref{rel1}, is implied by the vanishing of O'Neill's integrability tensor $A$. This is a rather strong condition, implying that $J^{-1}(\RR_+\mu)$ is a locally  (not necessarily Riemannian) product and cannot be predicted by the action. Other conditions on the second fundamental form which are common in Riemannian and Cauchy-Riemann submanifold theory, see e.g. \cite{be}, (mixed totally geodesic, (contact)-totally umbilical, extrinsic sphere etc.) and permit some speculations in \eqref{final} or even the computation of the Ricci curvature of the quotient, seem to be artificial in this context, as not directly expressible in terms of the action.

We apply the above computation for $N$ being $J^{-1}(0)$  and for
$P$ being the respective reduced space. Then  equation
\eqref{final} implies:
\begin{pr}
The reduced space at $0$ of a  Sasakian manifold with positive
$\varphi$-sectional curvature (in particular of an odd sphere with
the standard Sasakian structure) has strictly positive
$\varphi$-sectional curvature.
\end{pr}


\begin{thebibliography}{99}

\bibitem{am} R. Abraham, J. Marsden, \emph{Foundations of Mechanics},
 second edition, New York, Benjamin/Cummings, 1978.

\bibitem{al} C. Albert, \emph{Le th\'eor\`eme de r\'eduction de
 Marsden-Weinstein en g\'eom\'etrie
 cosymplectique et de contact}, {J. Geom. Physics}, {\bf 6} (1989),
627--649.

\bibitem{ar} V.I. Arnold, \emph{Mathematical Methods of Classical
 Mechanics},
 Springer Verlag, 1984.

\bibitem{be} A. Bejancu, \emph{The geometry of CR submanifolds},
 Mathematics
and its Applications (East European Series), 23. D. Reidel Publishing
 Co.,
Dordrecht, 1986.

\bibitem{bes} A. Besse, \emph{Einstein manifolds}, Springer Verlag
 (1987).

\bibitem {bl} D.E.  Blair, \emph{Riemannian geometry of contact and
symplectic manifolds}, Progress in Math. {\bf 203},
Birkh{\"a}user, Boston, Basel, 2002.

\bibitem{bg_se}Ch.P. Boyer, K. Galicki, \emph{On Sasakian-Einstein
 geometry}, Internat. J. Math.  {\bf 11}  (2000), 873--909.


\bibitem{dor} O. Dr\u agulete, L. Ornea, T.S. Ratiu, \emph{Cosphere
 bundle reduction in contact geometry}, J. Symplectic Geom. (to
 appear).

\bibitem{ge} H. Geiges, \emph{Constructions of contact manifolds},
 {Math. Proc. Cambridge Philos. Soc.}, {\bf 121} (1997), 455--464.

\bibitem{go} G. Grantcharov, L. Ornea, \emph{Reduction of Sasakian
 manifolds}, J. Math. Physics, {\bf 42} (2001), 3808--3816.

\bibitem{ko} Y. Kamishima, L. Ornea, \emph{Geometric flow on compact locally
conformal K\"ahler manifolds}, Tohoku Math. J. (to appear).

\bibitem{kob} S. Kobayashi,
\emph{Submersions of CR submanifolds},
Tohoku Math. J. {\bf 39} (1987), 95--100.

 \bibitem{mr1} J.E. Marsden, T.S. Ratiu, \emph{Introduction
 to Mechanics and Symmetry}, Springer Texts in Appl. Math. {\bf 17},
 Second edition, 1999.

  \bibitem{ov} L. Ornea, M. Verbitsky, \emph{Immersion theorem for
 compact
 Vaisman manifolds}, Math. Ann. (to appear), available at math.AG/0306077.

\bibitem{pap} N. Papaghiuc, \emph{Submersions of semi-invariant
 submanifolds of
a Sasakian manifold},  An. \c Stiin\c t. Univ. "Al. I. Cuza Ia\c
si", Sec\c t. I  Mat. {\bf 35} (1989),  281--288.

\bibitem{wil} C. Willett, \emph{Contact reduction},  Trans. Amer. Math.
 Soc.
{\bf 354}  (2002), 4245--4260.


 \end{thebibliography}
 \end{document}